\documentclass[11pt]{article}
\usepackage{enumerate}
\usepackage{amsmath,amsthm}
\usepackage{amssymb, latexsym}
\usepackage[mathscr]{euscript}
\usepackage{color}
\usepackage{array}
\usepackage{flafter}
\usepackage{layout}
\usepackage{graphicx}
\usepackage[applemac]{inputenc}
\usepackage{hyperref}
\theoremstyle{definition}
  \newtheorem{defi}{Definition}[section]
  
  \newtheorem{teo}{Theorem}[section]
  
  \newtheorem{lema}{Lemma}[section]
  \newtheorem{prop}{Proposition}[section]

\title{A note on large deviations for the stable marriage of Poisson and Lebesgue with random appetites}
\author{Daniel Andrés Díaz Pachón\thanks{Institute of Mathematical Statistics - Sao Paulo University. danielp@ime.usp.br}} 
\date{\today}
\allowdisplaybreaks
\begin{document}
\maketitle
\begin{abstract}
\noindent Let $\Xi\subset\mathbb R^d$ be a set of centers chosen according to a Poisson point process in $\mathbb R^d$. Let $\psi$ be an allocation of $\mathbb R^d$ to $\Xi$ in the sense of the Gale-Shapley marriage problem, with the additional feature that every center $\xi\in\Xi$ has an appetite given by a nonnegative random variable $\alpha$. Generalizing some previous results, we study large deviations for the distance of a typical point $x\in\mathbb R^d$ to its center $\psi(x)\in\Xi$, subject to some restrictions on the moments of $\alpha$.
\end{abstract}

\section{Introduction}

Let $\Xi\subset\mathbb R^d$ be a set of \textit{centers} chosen according to a Poisson point process of law $\textbf P$ and expectation $\textbf E$, whose intensity is $\lambda:=\textbf E\left[\Xi\cap[0,1]^d\right]$. Let $\alpha$ be a nonnegative random variable (r.v.) with law $\mathbb P$ and expectation $\mathbb E$ called \textit{appetite}. Every center $\xi\in\Xi$ will have a random appetite in such a way that appetites corresponding to different centers are independent of each other and identically distributed. The appetites are also independent of $\Xi$. With a little abuse of notation, but keeping in mind the independence, we will denote the appetite of a center $\xi$ by $\alpha_\xi$, when it does not lead to confusion (a formal construction is presented below). The elements $x\in\mathbb R^d$ will be called \textit{sites}.  Let $\psi$ be a function of $\mathbb R^d$ to $\Xi\cup\{\infty,\Delta\}$ called \textit{allocation}. The goal of $\psi$ is to allocate every site to a center. If $\psi$ cannot warrant that allocation, then the site will be sent to the element $\infty$. Also, there is a null set of sites under the Lebesgue measure $\mathcal L$ which will be sent to $\Delta$.

The \textit{territory} of the center $\xi$ will be given by $\psi^{-1}(\xi)$ and it satisfies the condition that $\mathcal L\left[\psi^{-1}(\xi)\right]\leq\alpha_\xi$. Whenever $\mathcal L\left[\psi^{-1}(\xi)\right]<\alpha_\xi$ the center will be \textit{unsated}; if $\mathcal L\left[\psi^{-1}(\xi)\right]=\alpha_\xi$ it will be \textit{sated}. Whenever $\psi(x)=\xi$ for some $\xi\notin\Delta$ we will say that $x$ is \textit{claimed}; if $\psi(x)=\infty$ we will call it \textit{unclaimed}.

Consider a site $x$ and a center $\xi$ such that $\psi(x)\notin\{\xi,\Delta\}$. We will say that the site $x$ \textit{desires} the center $\xi$ if $|x-\xi|<|x-\psi(x)|$ or if $\psi(x)=\infty$ (here and henceforth $|\cdot|$ is the Euclidean norm). If there exists $x'\in\psi^{-1}(\xi)$ such that $|x-\xi|<|x'-\xi|$, or if $\xi$ is unsated, we will say that the center $\xi$ \textit{covets} the site $x$. A pair $(x,\xi)$ will be \textit{unstable} if $x$ desires $\xi$ and $\xi$ covets $x$. If $\psi$ does not produce any unstable pair it will be called \textit{stable}. Finally, call $\mathcal C$ the closure of the set of claimed sites.

Hoffman, Holroyd, and Peres gave an explicit construction of the function $\psi$ for the case of constant appetite for every center in \cite{HHP1} and it was extended for the general case of random appetites in \cite{Diazphd}. Formally, it is done using the stable marriage algorithm of Gale-Shapley (see \cite{GS}). Informally, all the centers will grow a ball at the same linear speed, claiming for each one every site not belonging to the territory of any other center. The process continues for each center $\xi$ up to the point in which $\xi$ is sated (has claimed for itself a territory with volume $\alpha_\xi$). If its appetite is never satisfied the ball will continue growing to $\infty$ trying to find unclaimed sites and the center will be unsated.
\\

We proceed now to describe the probability space based on a construction proposed by Mester and Roy for the Poisson Boolean model (see \cite{MR}, pp. 16--17): Let $\Xi$ be defined on a probability space $(\Omega_1, \mathcal F_1, \textbf P)$. Let $(\Omega_2, \mathcal F_2, P)$ be the probability space on which the \textit{iid} r.v.'s $\{\alpha(n,z),n\in\mathbb N,z\in\mathbb Z^d\}$, distributed as $\alpha$, are defined (thus, $\Omega_2$ has a product structure such that each marginal has law $\mathbb P$).  Set $\Omega=\Omega_1\times\Omega_2$, equip it with the usual product $\sigma$-algebra and product measure $\mathcal P=\textbf P\times P$. The configuration of appetites of centers corresponding to $(\omega_1,\omega_2)$ is obtained using \textit{binary cubes of order} $n$:
\begin{align*}
	K(n,z):=\prod_{i=1}^d\left(z_i2^{-n},(z_i+1)2^{-n}\right] \text{ for all $n\in\mathbb N$ and $z\in				\mathbb Z^d$}.
\end{align*}

Now, each center $\xi\in\Xi$ is contained in a unique binary cube of order $n$, $K(n,z(n,\xi))$, and \textbf P-a.s., for each $\xi\in\Xi$ there is a unique smallest number $n_0=n_0(\xi)$ such that $K(n_0,z(n_0,\xi))$ contains no other centers of $\Xi$. The appetite of the ball centered at $\xi$ is given by $\alpha(n_0, z(n_0,\xi))$, which will be noted for simplicity as $\alpha_\xi$ when it does not lead to confusion, as it was said at the beginning.

Thus, the product structure of $\Omega$ warrants the independence between appetites and centers and the product structure of $\Omega_2$ warrants the independence among the appetites for different centers.

\textbf{Remark}. As it is noted in \cite{MR}, maybe the first construction that comes to mind is one in which the centers are ordered linearly following some determined law, and to construct in one probability space the point process $\Xi$ and \textit{iid} r.v.'s $\alpha_1,\alpha_2,\ldots$, and construct a realization by assigning $\alpha_i$ to the $i$-th center of $\Xi$. The problem with this setting is that we will need to use some ergodic results bur the ordering of the centers is not preserved with translations.
\\

The model can be denoted as $(\Xi,\alpha,\psi)$ with joint law $\mathcal P$ and expectation $\mathcal E$. Since $\psi$ is \textit{translation-equivariant} (i.e., $\psi_\Xi(x)=\xi$ implies that $\psi_{T\Xi}(Tx)=T\xi$, for every translation $T\in\mathbb R^d$ and where $\psi_\Xi$ is an allocation to the support $[\Xi]$ of $\Xi$), we have that $\mathcal P$ is translation invariant. Its \textit{Palm version} $(\Xi^*,\alpha^*,\psi^*)$, with law $\mathcal P^*$ and expectation operator $\mathcal E^*$, is conditioned to have a center on the origin. Note that since $\alpha$ is independent of $\Xi$, then $\alpha^*$ and $\Xi^*$ are also independent and $\alpha=\alpha^*$. 

By a result of Thorisson in \cite{Thor} (which was the beginning of the current research in allocations), $\Xi$ and $\Xi^*$ can be coupled so that almost surely one is a translation of the other.  Therefore, $\psi_{\Xi^*}$ is defined $\mathcal P^*$-a.s. Finally, remember that since $\Xi$ is a Poisson process, $\Xi^*$ is a Poisson process with an added point at the origin. (See a similar reasoning in \cite{HHP1} p. 1256 and \cite{HPPS} pp. 8--9). 

Defined this way, there are at least two quantities natural to consider in the model: $X':=|x-\psi(x)|$, which by translation-invariance happens to have the same law as $X=|\psi(0)|$; and $R_{\psi_{\Xi}}(\xi)=\text{ess sup}_{x\in\psi^{-1}(\xi)}|\xi-x|$: the radius of $\psi^{-1}(\xi)$. Thus, the radius of a typical center can be defined as $R^*=R_{\psi_{\Xi^*}}(0)$. 

Continuing with what was done in \cite{Diazpsmpl} with percolation, we generalize here some results of the stable marriage of Poisson and Lebesgue when the appetite is random.  The results presented in this article were first obtained in \cite{HHP2} for the case of constant appetite. We use the same methods to obtain large deviations for the case of random appetites, but of course there is a weakening: generally, when the appetites are constant the decay on the laws of $R^*$ and $X$ is exponential in $(R^*)^d$ and $X^d$, respectively; but now, with random  appetites, the decay is polynomial. We will keep the notation as consonant with \cite{HHP2} as possible.\\

Besides the existence of the stable allocation for a discrete set of centers, in \cite{Diazphd} there were some generalizations of previous results in \cite{HHP1} and \cite{HHP2}. Since most of proofs of these results remain almost unchanged (but in the random appetite case the results will hold $\mathcal P$-a.s.), here we mention some of them without proof and refer the reader directly to \cite{HHP1}. 

\begin{teo}[\textbf{Almost sure uniqueness}]\label{unique}
	Let $\Xi$ be a Poisson point process with law \textbf P and finite intensity $\lambda$, where each 		center $\xi$ chooses its appetite with law $\mathbb P$, then there exists a unique stable allocation 	for $\mathcal L$-a.e.\ $x\in\mathbb R^d$, $\mathcal P$-a.s.
\end{teo}

\begin{teo}[\textbf{Phase transitions}]\label{phases}
	Let $\Xi$ be a Poisson point process with law \textbf P and intensity $\lambda\in(0,\infty)$.
	\begin{itemize}
		\item If $\lambda\mathbb E\alpha<1$ (\textbf{subcritical}), then all the centers are sated but 				there exists an infinite volume of unclaimed sites, $\mathcal P$-a.s.
		\item If $\lambda\mathbb E\alpha=1$ (\textbf{critical}), then all the centers are sated and $				\mathcal L$-a.e.\ site is claimed, $\mathcal P$-a.s.
		\item If $\lambda\mathbb E\alpha>1$ (\textbf{supercritical}), then there are unsated centers 				but $	\mathcal L$-a.e.\ site is claimed, $\mathcal P$-a.s.
	\end{itemize}
\end{teo}

A set $A$ is an \textit{essential subset} of $B$ if $\mathcal L(B\setminus A)=0$. An allocation will be called \textit{canonical} if for every site $x$ such that there exists $\xi\in\Xi\cup\infty$ satisfying that $x$ has a neighborhood which is an essential subset of $\psi^{-1}(\xi)$, then $\psi(x)=\xi$; if there is no such $\xi$, then $\psi(x)=\Delta$. 

\begin{teo}[\textbf{Canonical allocations}]\label{teocanon}
	It is possible to construct a stable allocation to $\Xi$ which is also canonical and minimizes the $		\mathcal L$-null set $\Delta$.
\end{teo}

Sets of centers satisfying that they have a unique stable allocation $\mathcal L$-a.e.\ and a unique canonical allocation are called \textit{benign} sets. In particular, $\Xi$ is benign by Theorems \ref{unique} and \ref{teocanon}.\\

Now, for sets of centers $\Xi_1,\Xi_2,\ldots$ and $\Xi$, write $\Xi_n\Rightarrow\Xi$ if for any compact $K\subset\mathbb R^d$ there exists $N$ such that for $n>N$ we have $\Xi_n\cap K=\Xi\cap K$. For allocations $\psi_1, \psi_2,\ldots$ and $\psi$ write $\psi_n\rightarrow\psi$ a.e., if for $\mathcal L$-a.e.\ $x\in\mathbb R^d$ we have $\psi_n(x)\rightarrow\psi(x)$ in the one-point compactification $\mathbb R^d\cup\{\infty\}$. 

Exactly as it was done with $\{\alpha_\xi\}_{\xi\in\Xi}$, for $i\in\mathbb N$ we are going to construct a family of \textit{iid} random variables $\{\alpha_i(\xi)\}_{\xi\in\mathbb N}$ independent of $\alpha$, independent of every set of centers, and with the same distribution as $\alpha$ (Formally, we repeat the process above to obtain r.v.'s $\alpha_i(n_0,z(n_0,\xi))$, for $n_0\in\mathbb N$ and $z\in\mathbb Z^d$). Then, to the center $\xi\in\Xi_i$ corresponds the appetite $\alpha_i(\xi)$. We can define the random variables $\hat\alpha_i(\xi)$ as follows:

\begin{align*}
	\hat\alpha_i(\xi)=\begin{cases}
		\alpha_i(\xi) &\text{ for all } \xi\notin\Xi\\
		\alpha_\xi &\text{ otherwise}.
	\end{cases}
\end{align*}

Thus, as required, $\{\hat\alpha_i(\xi)\}_{\xi\in\mathbb N}$ is an \textit{iid} set of r.v.'s distributed as $\alpha$. We are going to use the following natural extension of Theorem 5 in \cite{HHP2}:

\begin{teo}[\textbf{Continuity}]\label{continuity}
	Let $\Xi_1,\Xi_2,...$ and $\Xi$ be benign sets of centers with appetites given by the $				\hat\alpha_i(\xi)$ for $\xi\in\Xi_i$ in $\Xi_i$ and $\alpha$ for all the centers in $\Xi$. Write $\psi_n=		\psi_{\Xi_n}$ $\psi=\psi_\Xi$ for their canonical allocations. If $	\Xi_n	\Rightarrow\Xi$, then $			\psi_n\rightarrow\psi$ a.e., $\mathbb P$-a.s.
\end{teo}

\noindent \textbf{Remark}. Theorems \ref{unique}, \ref{phases}, and \ref{teocanon} have more general statements in \cite{HHP1} (and \cite{Diazphd}). Theorem \ref{unique} was proved for every translation-invariant point process, and Theorems \ref{phases} and \ref{teocanon} were proved for every point process ergodic under translations. For our interest, the versions so stated will suffice. The proof of Theorem \ref{continuity} is essentially the same as the proof of Theorem 5 in \cite{HHP2} once a classical coupling is constructed via the r.v.'s $\hat\alpha_i(\xi)$. \\

There are some additional things needed for the proofs but they will be introduced in situ. For now, what we have suffices to state the results in this article. Finally, with a little abuse of notation to make the article more readable, we will use $\mathcal P$ and $\mathcal E$, except when it will be explicitly clear that we are dealing with $\Xi$ and $\alpha$.

\section{Statement of the results}

Our first result determines explicit bounds of the appetites under certain restrictions of $\alpha$. The next two theorems establish results for the upper bounds of $X$ and $R^*$ in the uncritical phases.

\begin{teo}\label{alfaextremo}
	Let $\Xi$ be a Poisson process in $\mathbb R^d$ with intensity 1 and consider a stable allocation 		with appetite $\alpha$.
	
	\begin{itemize}
		\item[(\textit i)] Let $\alpha$ be a r.v.\ with mean $\mu>2^d$, $d>1$ and variance $\sigma^2$. 			Then, given $\delta\in(0,d)$, we have 
			\begin{align*}
				\mathcal EX^{d-\delta}<\infty.
			\end{align*}
		\item[(\textit{ii})] Let $\alpha$ be a r.v.\ bounded from above by $M$. Assume that $\alpha$ 				has mean $\mu<2^{-d}$ and variance $\sigma^2$. Then, for every $\varepsilon\in(0,d/				2)$, we have that $\mathcal E^*e^{c(R^*)^\varepsilon}<\infty$, for some $c>0$.
	\end{itemize}
\end{teo}

\begin{teo}[\textbf{Supercritical upper bound}]\label{superc}
	Let $\Xi$ be a set of centers of intensity $\lambda>1$  and let $\alpha$ be the appetite with mean 1 	and variance $\sigma^2$. Then, for $\beta\in(0,d)$, we have that $\mathcal E X^{d-\beta}<\infty$.
\end{teo}

 \begin{teo}[\textbf{Subcritical upper bound}]\label{csubc}
 	Let $\Xi$ be a Poisson process of intensity $\lambda<1$ in $\mathbb R^d$ with $d>2$. Let $\alpha	$ be the appetite with mean 1 and finite variance $\sigma^2$. Then, for $\beta\in(0,d/2)$, we have 		that $\mathcal E^*(R^*)^{d/2-\beta}<	\infty$.
 \end{teo}

\section{Explicit bounds}

We start from mentioning two results of Nagaev for the large deviation of the sums of independent r.v.\ and two known results of Chernoff on the large deviations of a Poisson r.v.\\

\begin{prop}[\textbf{Corollary 1.8 in \cite{Nag}}]\label{cotaNagaev}
	Let $(X_i)_{i\in\mathbb N}$ be a sequence of \textit{iid} r.v.\ with $\mathcal EX_1=0$ and variance 		$\sigma^2$. Define $A_2^+:=\int_{u\geq0}u^2dF(u)$ and take $x>0$, then 
	\begin{equation*}
		\mathcal P\left[\sum_{i=1}^nX_i\geq x\right]\leq\frac{4nA_2^+}{x^2}+
		\exp\left\{\frac{-8x^2}{ne^2\sigma^2}\right\}
	\end{equation*}
	
\end{prop}
\begin{prop}[\textbf{Corollary 1.12 in \cite{Nag}}]\label{cotaNagaev2}
	Let $(X_i)_{i\in\mathbb N}$ be a sequence of \textit{iid} r.v.\ with $\mathcal EX_1=0$, $X_1\leq M$, 		and variance $\sigma^2$. Let $S_n=\sum_{i=1}^nX_i$ and take $x>0$, then
	\begin{equation*}
		\mathcal P[S_n\geq x]\leq\exp\left[x/M-(x/M+n\sigma^2/M^2)\log\left(\frac{xM}{n\sigma^2} 			+1\right)\right].
	\end{equation*}
\end{prop}

\begin{lema}[\textbf{Chernoff inequalities for Poisson r.v.}]\label{ChePo}
	Let $N$ be a Poisson r.v.  with mean $\lambda\in(0,+\infty)$, then
	\begin{align}
		\mathcal P[N\geq a]&\leq e^{-\lambda g(\lambda/a)}\text{\ \ \ for } a>\lambda\label{ChePo1}\\
		\mathcal P[N\leq a']&\leq e^{-\lambda g(\lambda/a')}\text{\ \ \ for } a'<\lambda\label{ChePo2}
	\end{align}
	
	\noindent where $g(x)=(x-1-\log x)/x$ (so $g(x)\rightarrow\infty$ when $x\rightarrow0$).
\end{lema}

\begin{proof}[Proof of Theorem \ref{alfaextremo}]
	To prove (\textit i) let $Z$ be the amount of centers of $\Xi$ in the ball $B(0,r)$, $r>0$. Thus, $Z$ 		is Poisson with mean $\pi_dr^d$, where $\pi_d$ is the volume of the unit ball in $d$ dimensions. 		We prove first  that
	
	\begin{align}\label{implic}
		\sum_{\xi\in B(0,r)}\alpha_\xi\geq\pi_d2^dr^d\Rightarrow X<r.
	\end{align}
	Suppose on the contrary that $\psi(0)=X\geq r$. Then 0 desires every center inside $B(0,r)$. 			Taking now the antecedent in (\ref{implic}) we have two options: Either every center in $B(0,r)$ is 		unsated, or there is at least one center $\xi\in B(0,r)$ with $\mathcal L$-non-null territory $T$ 			outside $B(0,2r)$; i.e., $T\subset\psi^{-1}(\xi)\cap[\mathbb R^d\setminus B(0,2r)]$ and $\mathcal 		LT>0$ for some $\xi\in B(0,r)$. In the first case, since every center $\xi'\in B(0,r)$ is unsated, then $		\xi'$ covets 0 by definition, thus $(0,\xi')$ is an unstable pair.  In the second case, for every $x\in T$ 		we have that $|\xi-x|>r>|\xi|=|\xi-0|$. Thus, $\xi$ covets 0 and $(0,\xi)$ is an unstable pair. A 			contradiction, and (\ref{implic}) is established. Then,
	
	\begin{align}\label{ineqset}
		\mathcal P[X>r]&\leq\mathcal P\left[\sum_{\xi\in B(0,r)}\alpha_\xi<\pi_d2^dr^d\right]\nonumber\\
		&\leq\mathcal P\left[\sum_{i=1}^Z\alpha_i\leq\pi_d2^dr^d\, \Big|\, Z\geq K\right]\mathcal 					P[Z\geq K]+\mathcal P[Z<K]\nonumber\\
		&\leq\mathcal P\left[\sum_{i=1}^{\lfloor K\rfloor}\alpha_i\leq\pi_d2^dr^d\right]+\mathcal P[Z<K]				\nonumber\\
		&=\mathcal P\left[\sum_{i=1}^{\lfloor K\rfloor}-\alpha_i'\geq-a\right]+\mathcal P[Z<K]\nonumber				\\
		&\leq\frac{4\lfloor K\rfloor A_2^+}{(-a)^2}+\exp\left\{\frac{-8(-a)^2}{\lfloor K\rfloor 						e^2\sigma^2}\right\}+e^{-\pi_dr^dg(\pi_dr^d/\lfloor K\rfloor)},
	\end{align}
	
	\noindent where $\alpha_k'=\alpha_k-\mu$ for $k=0,1,\ldots,K$, $K=\pi_d(r^d-r^{d/2+\varepsilon})$ 	for small positive $\varepsilon$;  $a=\pi_d2^dr^d-K\mu=\pi_dr^d(2^d-\mu-\mu r^{\varepsilon-d/		2})$. 	Then, by the defintion of $\mu$, $a<0$ when $r^{d/2-\varepsilon}>\mu/(\mu-2^d)$. The first 		two terms in (\ref{ineqset}) result from Proposition \ref{cotaNagaev} and the third term results from 		the Chernoff inequality (\ref{ChePo2}) in Lemma \ref{ChePo}. Given that $g(x)\sim(x-1)^2/2$ when 		$x\rightarrow1$, we have that $g\left(r^d/(r^d-r^{d/2+\varepsilon})\right)\sim r^{-d+2\varepsilon}/2$ 		when $r\rightarrow\infty$. 	Thus, for $r$ large enough, 
	
	\begin{align*}
		\mathcal P[X>r]&\leq c_1r^{-d}+e^{-c_2r^d}+e^{-c_3r^{2\varepsilon}}\\
		&\leq3c_1r^{-d}.\\
	\end{align*}
	
	\noindent where $c_3<\pi_d/2$, $c_2<8/(e^2\sigma^2)$, and $c_1<4A_2^+$. To finish the proof 		note 	that, given $\delta\in(0,d)$,
	\begin{align}\label{nicebound}
		\mathcal EX^{d-\delta} =\int_0^\infty\mathcal P\left[X\geq t^{1/(d-\delta)}\right]dt 						\leq1+3c_1\int_1^\infty t^{-d/(d-\delta)}dt<\infty.
	\end{align}
	
	In order to prove (\textit{ii}) let $Z'$ be the amount of centers of $\Xi^*$ in $B(0, 2r)$. Then 			$Z'-1$ is distributed as a Poisson with intensity $\pi_d2^dr^d$. If $Z'<M^{-1}\pi_dr^d$ we have 		that $\sum_{\xi\in B(0,2r)}\alpha_\xi<\pi_dr^d$. Therefore, there is a $\mathcal L$-non-null set $L		$ in $B(0,r)$ such that its elements do not belong to the territory of the centers of $\Xi^*$ in 			$B(0,2r)\cap\Xi^*$.  Hence, the stability of the pair $(x,0)$ implies that the appetite of 0 was 			already satisfied inside the closed ball $B[0,|x|]$. Then
	\begin{equation*}
		R^*<|x|<r.
	\end{equation*}
	
	\noindent We proceed as in part (\textit i), to obtain that
	\begin{align}\label{logar}
		\mathcal P[R^*>r]&\leq \mathcal P\left[\sum_{k=0}^{Z'}\alpha_k>\pi_dr^d\right]\nonumber \\ 					\nonumber
		&\leq\mathcal P\left[\sum_{k=0}^{Z'}\alpha_k>\pi_dr^d\,  \Big|\,  Z'\leq K\right]\mathcal P[Z'\leq 					K]+\mathcal P[Z'>K]\\
		&\leq\mathcal P\left[\sum_{k=0}^{\lceil K\rceil}\alpha_k'>a\right]+\mathcal P[Z'>K]
	\end{align}
	
	\noindent where $\alpha_k'=\alpha_k-\mu$ for $k=0,1,\ldots,K$; $K=\pi_d2^d(r^d+r^{d/2+			\varepsilon})$; $a=\pi_dr^d-K\mu=\pi_d r^d(1-2^d\mu(1+r^{-d/2+\varepsilon}))$; and $\varepsilon		\in(0,d/2)$. Then $a>0$ when $r^{d/2-\varepsilon}>2^d\mu/(1-2^d\mu)$. The appetite of the center 		at the origin is noted by $\alpha_0$. Now, after Proposition \ref{cotaNagaev2},the first term in 			(\ref{logar}) can be bounded by
	
	\begin{align}\label{logar1}
		\exp\left\{a/M-\left[a/M+\left(\lceil K\rceil+1\right)\sigma^2/M^2\right]\log\left(\frac{aM}						{\left(\lceil K\rceil+1\right)\sigma^2}+1\right)\right\},
	\end{align}
	
	\noindent where we will assume without loss of generality that $M>(e-1)\sigma^2\left(\left\lceil K		\right\rceil+1\right)/a$, in order to make the expression inside the logarithm greater than 1. Also, 		using Lemma \ref{ChePo}, the second term in (\ref{logar}) can be bounded by
	\begin{align}\label{logar2}
		\exp\{-\pi_d2^dr^d g(\pi_d2^dr^d/K)\}.
	\end{align}
	
	\noindent Note that $g(r^d/(r^d+r^{d/2+\varepsilon}))\sim r^{-d+2\varepsilon}/2$. With those 			considerations in mind, using (\ref{logar}), (\ref{logar1}), and (\ref{logar2}), we obtain
	\begin{align*}
		\mathcal P[R^*>r]&\leq e^{a/M-\left[a/M+K\sigma^2/M^2\right]}+e^{-c_1r^\varepsilon}\\
		&\leq e^{-c_2r^d}+e^{-c_1r^\varepsilon}\\
		&\leq 2e^{-c_3r^\varepsilon}.
	\end{align*}
	
	\noindent where $c_1=c_1(d)<\pi_d2^{d-1}$, $c_2=c_2(d)<\pi_d2^d\sigma^2/M^2$, and $c_3=			\min\{c_2,c_1\}$.\\
\end{proof}

\section{Supercritical bound}

Given a benign set $\Xi\subset\mathbb R^d$ and a measurable set $A\subset\mathbb R^d$, let $\Xi_A'$ be a random set of centers which is the union of $\Xi\cap A$ and a Poisson process with intensity $\lambda$ in $\mathbb R^d\setminus A$. Call $F_{\lambda,A}$ the law of $\Xi_A'$ and let $Q(L):=[-L,L)^d$. With these concepts at hand we will introduce the following definition and results which will equip us to prove Theorem \ref{superc}.

\begin{defi}[\textbf{Decisive set}]\label{deci}
 	Given a benign set $\Xi$, the measurable set $A\subset\mathbb R^d$ is $\Xi$-\textbf{decisive} for a 	site $x$ if for every $\lambda\in(0,\infty)$ we have that $F_{\lambda,A}$-a.s., $\Xi_A'$ is a benign 		set and $\psi_{\Xi_A'}(x)=\psi_\Xi(x)$.
\end{defi}

\begin{prop}\label{pdeci}
	Let $\hat\alpha_n(\xi)$ and $\alpha$, $\Xi_n$ and $\Xi$, and $\psi_n$ and $\psi$, be the appetites, 	sets of centers and allocations, respectively, as in Theorem \ref{continuity}. If $\psi(z)=\xi\in\Xi$ and 	$z$ is not equidistant from any two centers of $\Xi$, then $\psi_n(z)\rightarrow\xi$.
\end{prop}

\begin{proof}
	It corresponds to Proposition 19 in \cite{HHP2}. It is a refinement of Theorem \ref{continuity}.
\end{proof}

\begin{lema} \label{ldeci}
	Let $\mathbb E\alpha=1$ and let $\Xi$ be a Poisson process of intensity $\lambda\geq1$. Then, $		\textbf P$-a.s.\ and $\mathbb P$-a.s., there exists $L<\infty$ such that $Q(L)$ is $\Xi$-decisive for 0. 
\end{lema}

\begin{proof}
	Since $\lambda\geq1$ we have by Theorem \ref{phases} that 0 is claimed a.s.\ Also, 0 is not 			equidistant from any two centers, a.s. On the complement of the event in the lemma, for each $L$ 		there exists a benign $\Xi_L$ that agrees with $\Xi$ on $Q(L)$ such that $\psi_{\Xi_L}(0)\neq\psi_		\Xi(0)$. But by Proposition \ref{pdeci} this is an event with null probability $\mathbb P$, since when 	this happens we have $\psi_{\Xi_L}(0)\rightarrow\psi_\Xi(0)$ as $L\rightarrow\infty$, $\mathbb P$-		a.s.
\end{proof}

\begin{lema}[\textbf{Decisive cubes}]\label{qdeci}
	Let $\mathbb E\alpha=1$ and let $\Xi$ be a Poisson process of intensity $\lambda\geq1$. For 		every $\varepsilon>0$ there exists $M<\infty$ such that
	\begin{align*}
		\mathcal E\mathcal L\left[x\in Q(M):Q(M) \text{ is not } [\Xi]\text{-decisive for } x\right]<				\varepsilon(2M)^d.
	\end{align*}
\end{lema}

\begin{proof}
	Fix $\varepsilon>0$. Let $U^L$ be the random set of sites $x$ for which $Q(L)+x$ is not $\Xi$-decisive. Then $U^L$ is translation-invariant in law, and by Lemma \ref{ldeci}, we may fix $L$ large 	enough so that it has intensity less than $\varepsilon/2$. For $M$ sufficiently large we have
	\begin{equation*}
		\textbf E\mathcal L\{U^L\cap Q(M)\}<(\varepsilon/2)(2M)^d,
	\end{equation*}
	thus
	\begin{align*}
		\textbf E\mathcal L&\{x\in Q(M):Q(M) \text{ is not } \Xi\text{-decisive for } x\}\\
		&<(\varepsilon/2)(2M)^d+(2M)^d-(2M-2L)^d\\
		&<\varepsilon(2M)^d.
	\end{align*}
\end{proof}

Lemmas \ref{ldeci}, \ref{qdeci} and Proposition \ref{pdeci} were introduced in \cite{HHP2}. We transcribed them here for the sake of thoroughness and contextualization. Now we have what we need to proceed with the proof of Theorem \ref{superc}:
 
\begin{proof}[Proof of Theorem \ref{superc}]
	Note first that if
	\begin{align}\label{unsated}
		\text{there exists } \xi\in[\Xi]\cap B(0,r) \text{ with } \mathcal L\left[\psi^{-1}(\xi)\cap B(0,2r)\right]<			\alpha_\xi,
	\end{align}
	
	\noindent then $X\leq r$. Otherwise, if $X>r$, then $\xi$ covets 0 and 0 desires $\xi$, and thus 			$(0,\xi)$ is unstable. Hence, 
	\begin{align*}
		\mathcal P[X>r]\leq\mathcal P\left\{\mathcal L\left[\psi^{-1}(\xi)\cap B(0,2r)\right]=\alpha_\xi,			\text{ for all } \xi\in\Xi\cap B(0,r)\right\}.
	\end{align*}
	
	Thus, it is enough to prove that the event at the right hand of the former inequality has polynomial 		decay. For a given $\lambda$ let
	\begin{align*}
		\varepsilon=\frac{\lambda-1}{10\cdot2^d}\wedge1,
	\end{align*}
	
	\noindent and let $M=M(\lambda,\varepsilon)$, as in Lemma \ref{qdeci}.
	
	For all $r>0$ we are going to place disjoint cubes congruent to $Q(M)$ in $B(0,2r)\setminus B(0,r)		$. Define $Q_z=Q(M)+2Mz$, for $z\in\mathbb Z^d$  and define
	\begin{equation*}
		Y_z=\mathcal L(x\in Q_z: Q_z \text{ is not } \Xi\text{-decisive for } x).
	\end{equation*}
	
	Define also
	\begin{align*}
		I=I(r)=\left\{z\in\mathbb Z^d:Q_z\subset B(0,2r)\setminus B(0,r)\right\},
	\end{align*}
	
	\noindent this is the index set of the cubes which are completely inside the shell. Finally, let
	\begin{align*}
		S=S(r)=\left[B(0,2r)\setminus B(0,r)\right]\setminus\bigcup_{z\in I}Q_z
	\end{align*}
	
	\noindent be the part of the shell which is uncovered by $\bigcup_{z\in I}Q_z$. For large enough $r		$ we have that
	\begin{align}\label{LS}
		\mathcal LS<\varepsilon\pi_dr^d.
	\end{align}
	
	\noindent To prove this, it suffices to prove that 
	\begin{align*}
		\mathcal L\left[B(0,r+\sqrt2M) \setminus B(0,r)\right]+\mathcal L\left[B(0,2r)\setminus B(0,2r-			\sqrt2M)\right]<\varepsilon\pi_dr^d,
	\end{align*} 
	for large enough $r$. Which is equivalent to prove that $\pi_d[(r+\sqrt2M)^d-r^d]+\pi_d[2^dr^d-(2r-		\sqrt2M)^d]<\varepsilon\pi_dr^d$ for large enough $r$. But this is clear since every factor of $\pi_d$ 	at the left side of the inequality has order $r^{d-1}$.
	
	Consider now the next events: 
	\begin{align*}
		E&=\left\{\sum_{\xi\in B(0,r)}\alpha_\xi>(\lambda-\varepsilon)\pi_dr^d\right\};\\
		G&=\left\{\sum_{z\in I}Y_z<4\varepsilon2^d\pi_dr^d\right\}
	\end{align*}
	and suppose that $E$ and $G$ occur. Then
		\begin{align*}
		\mathcal L[x\in B(0,2r):\psi(x)\in B(0,r)]&\leq\sum_{z\in I}Y_z+\pi_dr^d+\mathcal LS\\
		&\leq(4\varepsilon2^d+1+\varepsilon)\pi_dr^d\\
		&<(\lambda-\varepsilon)\pi_dr^d\\
		&<\sum_{\xi\in B(0,r)}\alpha_\xi,
	\end{align*}
	
	\noindent where the third inequality is obtained because of the definition of $\varepsilon$. Then 		the event in (\ref{unsated}) occurs when $E$ and $G$ occur. So it is enough to study the tails of $		\mathcal P[E^C]$ and $\mathcal P[G^C]$ to conclude that the decay of $X$ will be at least the one 		of these two with heavier tail. 
	
	We start with $\mathcal P[G^C]$. Note that $(Y_z)_{z\in I}$ is a set of \textit{iid} r.v. Furthermore, 		$Y_z\in(0,2^dM^d)$. Thus, $Y_z$ has finite mean $\mu_Y$ (in fact, $\mu_Y<\varepsilon(2M)^d/2$) 	and variance $\sigma^2_Y$. Let $\#I$ be the cardinality of the set $I$ (note that $\#I$ has order 		$r^d$) and let $K=\pi_dr^d(2^d-1)/(2M)^d$, which would be the maximum number of cubes with 		side $2M$ in the shell $B(0,2r)\setminus B(0,r)$. Then, using the Nagaev inequality of Proposition 		\ref{cotaNagaev}, we obtain that, for large $r$,
	\begin{align*}
		\mathcal P\left[\sum_{z\in I}Y'_z\geq y\right]&\leq\mathcal P\left[\sum_{z=1}^{\lceil K\rceil}Y'_z					\geq y\right]\\
		&\leq4\lceil K\rceil A_2^+/y^2+\exp\{-8y^2/(\lceil K\rceil e^2\sigma^2_Y)\}\\
		&\leq c'r^{-d}
	\end{align*}
	
	\noindent where $Y'_z=Y_z-\mu_Y$, $y=4\varepsilon2^d\pi_dr^d-(\#I)\mu_Y$.\\
	
	As for $\mathcal P[E^C]$, let $\alpha_k'=\alpha_k-1$ for $k=1,\ldots,K$, $K=\pi_d(\lambda r^d-r^{d/		2+\delta})$, where $\delta$ is a small enough positive real. Let $a=(\lambda-\varepsilon)\pi_dr^d-		K=\pi_dr^{d/2+\delta}(1-\varepsilon r^{d/2-\delta})	$, which will be negative for $r>\varepsilon^{-2/	(d-2\delta)}$. Finally, let $N$ be a Poisson r.v.\ with 	intensity $\lambda\pi_dr^d$. Then we have, 
	\begin{align*}
		\mathcal P\left[E^C\right]&=\mathcal P\left[\sum_{\xi\in B(0,r)}\alpha_\xi\leq(\lambda-					\varepsilon)\pi_dr^d\right]\\
		&\leq\mathcal P\left[\sum_{i=1}^N\alpha_i\leq(\lambda-\varepsilon)\pi_dr^d|N\geq K\right]					\mathcal P[N\geq K]+\mathcal P[N<K]\\
		&\leq\mathcal P\left[\sum_{i=1}^{\lfloor K\rfloor}\alpha_i\leq(\lambda-\varepsilon)\pi_dr^d						\right]\mathcal P[N\geq K]+\mathcal P[N<K]\\
		&\leq\mathcal P\left[\sum_{i=1}^{\lfloor K\rfloor}-\alpha_i'\geq -a\right]+\mathcal P[N<K].
	\end{align*}
	
	\noindent By the Chernoff inequality (\ref{ChePo2}) in Lemma \ref{ChePo} we obtain for the second 	term in the last inequality that 
	\begin{align}\label{gdChePo1}
		\mathcal P[N<K]&\leq e^{-\lambda\pi_dr^dg(\lambda\pi_dr^d/K)}\nonumber\\
		&\leq e^{-c_1r^{2\delta}}
	\end{align}
	
	\noindent for some $c_1>0$. And for the first term, using Proposition \ref{cotaNagaev}, we obtain 		that
	\begin{align}\label{gdNag1}
		\mathcal P\left[\sum_{i=1}^{\lfloor K\rfloor}-\alpha_i'\geq -a\right]&\leq\frac{4\lfloor K\rfloor A_2^			+}{a^2}+\exp\left\{\frac{-8a^2}{\lfloor K\rfloor e^2\sigma^2} \right\}\nonumber\\
		&\leq c_3r^{-d}+e^{-c_2r^d}.
	\end{align}
	
	\noindent Thus, from (\ref{gdChePo1}) and (\ref{gdNag1}) we obtain
	\begin{align*}
		\mathcal P\left[E^C\right]&\leq c_3r^{-d}+e^{-c_2r^d}+e^{-c_1r^d}\\
		&\leq c_4r^{-d}.
	\end{align*}
	
	\noindent Since both terms have polynomial decay, we obtain the result after an application of 		(\ref{nicebound}).
\end{proof}
 
 \section{Subcritical bound}
 
In order to prove Theorem \ref{csubc}, we will need Proposition \ref{recu} below. Although the proof of this Proposition is almost identical to the one stated by Hoffman, Holroyd and Peres (see Lemmas 14 and 18, and Corollary 15 in \cite{HHP2}), we reproduce it here for the sake of clarity and completeness, with the suitable minor changes needed in this case. 

\begin{lema}\label{sated}
	Suppose $\Xi_n\Rightarrow\Xi$ and let $\alpha'_n(\xi)$ be the appetite of the centers $\xi\in\Xi_n$ 		and $\alpha$ in $\Xi$ as in Theorem \ref{continuity}. Then $\psi_n\rightarrow\psi$, almost 			everywhere. If there is a set $A$ of positive volume such that every $z\in A$ desires $\xi$ under 		$\psi$, then for $n$ sufficiently large, $\xi$ is sated in $\psi_n$ and
	\begin{equation}\label{sated1}
		\limsup_{n\rightarrow\infty} R_{\psi_n}(\xi)\leq\text{ess\ inf}_{z\in A}|z-\xi|<\infty.
	\end{equation} 
\end{lema}

\begin{proof}
	It is identical to the proof of Lemma 18 in \cite{HHP2}.
\end{proof}

\begin{defi}[\textbf{Replete set}]\label{cheio}
 	Given a benign set $\Xi$, we say that a measurable set $A\subset\mathbb R^d$ is 					$\Xi$-\textbf{replete} for a center $\xi$ if for every $\lambda\in(0,\infty)$ we have that,  $				F_{\lambda,A}$-a.s., $\Xi_A'$ is benign and $\mathcal L\left(\psi_{\Xi_A'}^{-1}(\xi)\cap A				\right)=\alpha_\xi$.
\end{defi}

\begin{lema}\label{replete} 
	Let $\mathbb E\alpha=1$ and let $\Xi$ be a Poisson process of intensity $\lambda<1$. Call $G$ 		the event that for every $\xi\in\Xi$ there exists $L<\infty$ such that $\xi+Q(L)$ is $\Xi$-replete for $		\xi$. Then $\textbf P(G)=1$.
\end{lema}

\begin{proof}
	On $G^C$ there exists a center $\xi\in\Xi$ such that for each $L$ there exists a set of centers $		\Xi_L$ agreeing with $\Xi$ on $\xi+Q(L)$ satisfying
	\begin{equation}\label{repleteproof}
		\mathcal L[\psi_{\xi_L}^{-1}(\xi)\cap(\xi+Q(L))]<\alpha_\xi.
	\end{equation}
	Define the appetite for the set of centers $\Xi_L$ as in Theorem \ref{continuity}: by $\hat				\alpha_n(\xi)$ for all $\xi\in\Xi_L$. Since $\Xi_L\Rightarrow\Xi$, we have that $\psi_L\rightarrow		\psi_\Xi$, $\mathbb P$-as., by Theorem \ref{continuity}. Furthermore, Lemma \ref{sated} applies to 		$\xi$ (by the subcritical phase in Theorem \ref{phases}), 	so almost 	surely for $L$ large enough, $		\xi$ is sated in each $\psi_L$ and the radii $R_{\psi_L}(\xi)$ are bounded as $L\rightarrow\infty$. 		Thus $G^C$ is null a.s.
\end{proof}

\begin{prop}\label{recu}
	Let $\mathbb E\alpha=1$ and let $\Xi$ be a Poisson process with intensity $\lambda<1$. For all 		$\varepsilon>0$ there exists $M$ such that
	\begin{align*}
		\mathcal E\left\{\xi\in\Xi\cap Q(M): Q(M)\text{ is not } [\Xi]\text{-replete for } \xi\right\}<						\varepsilon(2M)^d.
	\end{align*} 
\end{prop}

\begin{proof}
	For $A\subset\mathbb R^d$, call $\Xi^L(A)$ the number of $\xi\in\Xi\cap A$ such that $\xi+Q(L)$ is 		not $\Xi$-replete for $\xi$. Lemma \ref{replete} and the monotone convergence theorem imply that 		$\textbf E\Xi^L(Q(1))\rightarrow0$ when $L\rightarrow\infty$. Thus we can choose $L<\infty$ so that 	the translation invariant point process $\Xi^L$ has intensity less than $\varepsilon/2$. Note that for 		$M>L$ and $\xi\in\Xi\cap Q(M-L)$, if $Q(M)$ is not $\Xi$-replete for $\xi$, then $\xi\in\Xi^L$. 			Therefore
	\begin{align*}
		\textbf E\#&\{\xi\in\Xi\cap Q(M):Q(M) \text{ is not } \Xi-\text{replete for } \xi\}\\
		&(\varepsilon/2)(2M-2L)^d+(2M)^d-(2M-2L)^d;
	\end{align*}
	which is smaller than $\varepsilon(2M)^d$ for large enough $M$.
\end{proof}

\begin{teo}\label{csubc2}
	Assume $\alpha$ has mean 1 and variance $\sigma^2$ and let $\Xi$ be a Poisson process with 		intensity $	\lambda<1$. Then there exists $c>0$ such that, for all $r>0$,
	\begin{align*}
		\mathcal P[\text{there exists } \xi\in\Xi\cap B(0,1)\text{ such that } R(\xi)>r]<cr^{-d}.
	\end{align*}
\end{teo}

\begin{proof}
	Fix $\lambda<1$ and note that if
	\begin{align}\label{R}
		\text{there exists } y\in B(0,r) \text{ with }\psi(y)\notin B(0,2r+1)
	\end{align}
	\noindent then $R(\xi)<r+1$ for all $\xi\in\Xi\cap B(0,1)$. Otherwise, $|y-\xi|<r+1$ and $|y-\psi(y)|>r		+1$ and $(y,\xi)$ would be unstable. Then it suffices to study the tail of the complementary event in 	(\ref{R}). Let $\varepsilon=(1-\lambda)/(10^d2^d)$ and let $M=M(\lambda,\varepsilon)$ as in 			Lemma \ref{recu}.
	
	Now, for all $r>0$ we are going to place disjoint copies of $Q(M)$ in the shell $B(0,2r+1)\setminus 		B(0,r)$. For $z\in\mathbb Z^d$ let $Q_z=Q(M)+2Mz$ and define
	\begin{equation*}
		W_z=\{\xi\in\Xi\cap Q_z:Q_z \text{ is not } \Xi\text{-replete for }\xi\}.
	\end{equation*}
	
	\noindent Let 
	\begin{align*}
		I=I(r)=\{z\in\mathbb Z^d:Q_z\subset B(0,2r+1)\setminus B(0,r)\}
	\end{align*}
	
	\noindent be the index set of the cubes completely contained  in the shell and let 
	\begin{align}
		S=S(r)=[B(0,2r+1)\setminus B(0,r)]\setminus\bigcup_{z\in I}Q_z
	\end{align}
	
	\noindent be the part of the shell which is uncovered by the cubes $Q_z$ with $z\in I$. Consider 		the following events:
	\begin{align*}
		E&=\left\{\sum_{\xi\in B(0,r)}\alpha_\xi<(\lambda+\varepsilon)\pi_dr^d\right\}\\
		F&=\left\{\sum_{\xi\in S}\alpha_\xi<\varepsilon\pi_dr^d\right\}\\
		G&=\left\{\sum_{z\in I}\sum_{\xi\in W_z}\alpha_\xi<4^d\varepsilon2^d\pi_dr^d\right\}.
	\end{align*}
	
	\noindent If $E$ occurs, we have that
	\begin{align*}
		\mathcal L\{y\in B(0,r):\psi(y)\notin B(0,r)\}&\geq\pi_dr^d-\sum_{\xi\in B(0,r)}\alpha_\xi\\
		&\geq\pi_dr^d-(\lambda+\varepsilon)\pi_dr^d\\
		&\geq9^d\varepsilon2^d\pi_dr^d,
	\end{align*}
	
	\noindent where the third inequality is due to the $\varepsilon$ chosen. If $F$ occurs, we have that
	\begin{align*}
		\mathcal L\{y\in B(0,r):\psi(y)\in S\}&<\sum_{\xi\in S}\alpha_\xi<\varepsilon\pi_dr^d.
	\end{align*}
	
	\noindent If $G$ occurs, by Definition \ref{cheio}, we have
	\begin{align*}
		\mathcal L\left\{y\in B(0,r):\psi(y)\in\bigcup_{z\in I}Q_z\right\}\leq\sum_{z\in I}\sum_{\xi\in W_z}				\alpha_\xi<4^d\varepsilon2^d\pi_dr^d.
	\end{align*}
	
	\noindent Since $B(0,2r+1)=S\cup B(0,r)\cup\bigcup_{z\in I}Q_z$, we have that if the event $E\cap 		F\cap G$ occurs, then
	\begin{align*}
		\mathcal L\{y\in B(0,r):\psi(y)\notin B(0, 2r+1)\}\geq(9^d2^d-4^d2^d-1)\varepsilon\pi_dr^d>0.
	\end{align*}
	
	\noindent Therefore, when $E$, $F$ and $G$ occur, the event in (\ref{R}) occurs too. So it is 			enough to analyze the decay of $\mathcal P[E^C]$, $\mathcal P[F^C]$ and $\mathcal P[G^C]$. We 	start with $G^C$. 
	
	Note that $\#I$ is bounded by $\pi_d\left[(2r+1)^d-r^d\right]/(2M)^d<\pi_d3^dr^d/(2M)^d$. Also, $		\mathcal E(\#W_z)$ is bounded by $\varepsilon(2M)^d$ because of Propositiion \ref{recu}. With 		this considerations, we define a Poisson r.v.\  $N$ having intensity $\pi_d3^dr^d\varepsilon$. 			Define $K=\pi_d(3^dr^d\varepsilon+r^{d/2+\delta})$, for some small $\delta>0$. Call $\alpha'_k=		\alpha_k-1$ for $k=1,\ldots,K$ 	and $a=4^d\varepsilon2^d\pi_dr^d-\lceil K\rceil$.
	\begin{align*}\label{cardI}
		\mathcal P[G^C]&=\mathcal P\left[\sum_{z\in I}\sum_{\xi\in W_z}\alpha_\xi\geq4^						\varepsilon2^d\pi_dr^d\right]\nonumber\\
		&\leq\mathcal P\left[\sum_{i=1}^N\alpha_i\geq4^d\varepsilon2^d\pi_dr^d\ |N\leq K\right]						\mathcal P[N\leq K]+\mathcal P[N>K]\\
		&\leq\mathcal P\left[\sum_{i=1}^{\lceil K\rceil}\alpha'_i\geq a\right]+\mathcal P[N>K]
	\end{align*}
	
	\noindent The second term in the last inequality can be bounded by $e^{-c_3r^{2\delta}}$ for some 	constant $c_3=c_3(\lambda,d)>0$, and the first term can be bounded using the Nagaev inequality 		in Proposition \ref{cotaNagaev}:
	\begin{align*}
		\mathcal P\left[\sum_{i=1}^{\lceil K\rceil}\alpha'_i>a\right]&\leq\frac{4\lceil K\rceil A_2^+}{a^2}+				\exp\left\{\frac{-8a^2}{\lceil K\rceil e^2\sigma^2}\right\}\\
		&\leq c_1r^{-d}+e^{-c_2r^d}.
	\end{align*}

	Now we consider $\mathcal P[E^C]$. Let $N'$ be a Poisson r.v.\ with mean $\lambda\pi_dr^d$. 		Call $\alpha_i'=\alpha_i-1$; $a=(\lambda+\varepsilon)\pi_dr^d-K$ where $K=\lambda\pi_d(r^d		+r^{d/2+\delta})$ for some small $\delta>0$. Note that $a$ is positive for $r>(\varepsilon\pi_d)^{-1/		d}$. Then
	\begin{align*}
		\mathcal P[E^C]&=\mathcal P\left[\sum_{\xi\in B(0,r)}\alpha_\xi\geq(\lambda+\varepsilon)\pi_d 			r^d\right]\\
		&\leq\mathcal P\left[\sum_{i=1}^{N'}\alpha_i\geq(\lambda+\varepsilon)\pi_dr^d|N'\leq K\right]				\mathcal P[N'\leq K]+\mathcal P[N'>K]\\
		&\leq\mathcal P\left[\sum_{i=1}^K\alpha_i'\geq a\right]+\mathcal P[N'>K].
	\end{align*}
	
	\noindent By the Chernoff inequality (\ref{ChePo1}), the second term is bounded by $e^{-			C_1r^{2\delta}}$, where $C_1=C_1(\lambda,d)>0$. By Proposition \ref{cotaNagaev}, the first term 		can be bounded as 	follows:
	\begin{align*}
		\mathcal P\left[\sum_{i=1}^K\alpha_i'\geq a\right]&\leq\frac{4\lceil K\rceil A_2^+}{a^2}+					\exp\left\{\frac{-8a^2}{\lceil K\rceil e^2\sigma^2}\right\}\\
		&\leq C_2r^{-d}+e^{-C_3r^d}.
	\end{align*}
	
	\noindent for some positive constants $C_1$ and $C_2$. Thus, $\mathcal P[E^C]$ and $\mathcal 		P[G^C]$ decay 	as $r^{-d}$. A similar procedure will show that $\mathcal P[F^C]$ has the same 		polynomial decay, which proves the theorem.
\end{proof}

\begin{proof}[Proof of Theorem \ref{csubc}]
	Let $Y$ be the number of centers $\xi$ in $\Xi\cap B(0,1)$ such that $R(\xi)>r$. Let $u=r^{d/2}$. 		Note that
	\begin{align*}
		\mathcal EY&=\mathcal EYI\{0<Y\leq u\}+\mathcal EYI\{Y>u\}\\
		&\leq u\mathcal P[Y>0]+\mathcal E[\Pi\cap B(0,1)]I\{\Pi\cap B(0,1)>r\}\\
		&\leq Cr^{-d/2}+\mathcal E[(\Pi\cap B(0,1))^2]/u\\
		&\leq Cr^{-d/2}+C_1r^{-d/2}
	\end{align*}
	
	\noindent where the second inequality is obtained because, by Theorem \ref{csubc2}, $\mathcal 		P[Y>0]$ is bounded by $Cr^{-d}$, for some $C>0$. The second term is bounded by $\mathcal E[\Xi		\cap B(0,1))^2]/u$. By a property of the Palm process we have that $\mathcal EY=\lambda\pi_d		\mathcal P^*[R^*>r]$. Finally, after applying (\ref{nicebound}), we obtain the result.
\end{proof}

\section{Acknowldegments}

The author is grateful to Serguei Popov, his Ph.D.\ advisor, for his guidance in the development of this work. He thanks the anonymous referee for the very useful suggestions and the commentator “biosdev” on a post of the author's blog (http://wp.me/sm7GZ-513). He also thanks Óscar Ocampo for some valuable insights. Finally, he acknowledges the financial support of CAPES during the development of his research.

\end{document}